\documentclass[10pt]{article}
\usepackage{amssymb, amsmath, amsfonts,latexsym, oldgerm,amscd,amsthm}
\usepackage{relsize}
\usepackage{url}

\textwidth13cm \usepackage[all]{xy} 
\newcommand{\tn}{\textnormal} \newcommand{\lra}{\longrightarrow}
 \newtheorem{de}{Definition}[section]

\newtheorem{re}[de]{Remark} 
\newtheorem{thm}[de]{Theorem} \newtheorem{lem}[de]{Lemma}
 \newtheorem{nt}[de]{Notation}
 \newtheorem{cor}[de]{Corollary}
 
\newcommand{\gm}{\mathfrak{m}}

\def\E{{\rm E}}

\def\ESp{{\rm ESp}}
\def\Sp{{\rm Sp}}

\def\GL{{\rm GL}}

\def\SL{{\rm SL}}

\def\lra{\longrightarrow}

\title{ Gram--Schmidt--Vaserstein generators \\ for odd sized
elementary groups} \author{Pratyusha Chattopadhyay 
~\& Ravi A. Rao 
\footnote{Correspondence author: {\it email:
ravi@math.tifr.res.in}} \\ {\small Statistics and Mathematics Unit,
Indian Statistical Institute, Bangalore 560 059, India} {\small \&}\\
{\small Tata Institute of Fundamental Research, Dr. Homi Bhabha Road,
Mumbai 400 005, India} }  \date{}

\begin{document}
\maketitle

\begin{center}{\it 2010 Mathematics Subject Classification: {13A99,
14L35, 15A24}}
\end{center}

\begin{center}{\it Key words: Elementary and symplectic matrices,
Local-Global principle of Quillen}
\end{center}

{\small ~~~~Abstract: The Gram-Schmidt process of L.N. Vaserstein for
making an elementary matrix symplectic yields a nice set of generators
for the odd sized elementary group.}

\section{\large Introduction}

The Gram--Schmidt process is a method for orthonormalising a set of
vectors in an inner product space. This enables one to transform a
square matrix (over the reals) by an upper triangular matrix to an
orthogonal matrix. 

A similar process exists in \cite{SV}, due to L.N. Vaserstein, in the
symplectic group $\Sp_{2n}(R)$. Here he showed that given an
elementary matrix (over a commutative ring $R$) of even size $2n$ one
can transform it by an elementary matrix of size $2n - 1$ to a
(elementary) symplectic matrix. 

Let $\varphi$ be an invertible alternating matrix of size $2n$. Then
L.N. Vaserstein's method actually permits one to transform an
elementary matrix of even size $2n$ by an elementary matrix of size
$2n - 1$ so that it is a (elementary) symplectic matrix
w.r.t. $\varphi$. 

We decided to collect all these odd sized elementary matrices which
transform a given elementary matrix to a symplectic matrix
w.r.t. $\varphi$. We call the subgroup of the elementary group
$\E_{2n-1}(R)$ generated by these as the group of elementary matrices
$\E_\varphi(R)$ w.r.t. $\varphi$. 

\vskip0.15in Our main result is the following observation:

\vskip0.15in

\noindent {\bf Theorem:} For an invertible alternating matrix
$\varphi$ of size $2n$, the subgroup $\E_\varphi(R)$ coincides with
$\E_{2n - 1}(R)$. 

\vskip0.15in This is proved by using a well-known Local-Global
argument of D. Quillen in \cite{Q}; to reduce the problem to the case
when $\varphi = \psi_r$, the standard symplectic matrix.

\vskip0.15in

We have also proved a relative version of the above theorem w.r.t. an
ideal of the ring $R$.

\section{\large Preliminaries}

\medskip A row ${ v} = (v_{1}, \ldots, v_{n}) \in R^{n}$ is said to be
{\it unimodular} if there are elements $w_{1}, \ldots, w_{n}$ in $R$
such that $v_{1}w_{1} + \cdots + v_{n}w_{n} = 1$. Um$_{n}(R)$ will
denote the set of all unimodular rows ${ v} \in R^{n}$. Let $I$ be an
ideal in $R$. We denote by ${\rm Um}_n(R, I)$ the set of all
unimodular rows of length $n$ which are congruent to $e_1 = (1, 0,
\ldots, 0)$ modulo $I$. (If $I = R$, then ${\rm Um}_n(R, I)$ is ${\rm
Um}_n(R)$).

Let E$_{n}(R)$ denote the subgroup of SL$_{n}(R)$ consisting of all
{\it elementary} matrices, i.e. those matrices which are a finite
product of the {\it elementary generators} E$_{ij}(\lambda) = I_{n} +
e_{ij}(\lambda)$, $1 \leq i \neq j \leq n$, $\lambda \in R$, where
$e_{ij}(\lambda) \in$ M$_{n}(R)$ has an entry $\lambda$ in its $(i,
j)$-th position.

In the sequel, if $\alpha$ denotes an $m \times n$ matrix, then we let
$\alpha^t$ denote its {\it transpose} matrix. This is of course an
$n\times m$ matrix. However, we will mostly be working with square
matrices, or rows and columns.

\begin{de} {\bf The Relative Groups E$_n(I)$, E$_n(R, I)$:}  {\rm Let
$I$ be an ideal of $R$. The {\it relative elementary group} {\rm
E}$_n(I)$ is the subgroup of {\rm E}$_n(R)$ generated as a group by
the elements {\rm E}$_{ij}(x)$, $x \in I$, $1 \leq i \neq j \leq n$.

The {\it relative elementary group} ${\rm E}_n(R, I)$ is the normal
closure of {\rm E}$_n(I)$ in {\rm E}$_n(R)$. }

\end{de}

\begin{lem} \label{suslin} {\rm (See \cite{18}, Corollary 1.2, Lemma
1.3)} Let $v, w \in R^n$, with $n \ge 3$ and $\langle v, w \rangle = v
\cdot w^t = 0$. Assume that $v$ is unimodular and $w \in I^{2n-1}
(\subseteq R^{2n-1})$. Then $I_n + v^t w \in E_n(R, I)$.
\end{lem}

\begin{nt} {\rm Let $M$ be a finitely presented $R$-module and $a$ be
a non-nilpotent element of $R$. Let $R_a$ denote the ring $R$
localized at the multiplicative set $\{a^i : i \ge 0 \}$ and $M_a$
denote the $R_a$-module $M$ localized at $\{a^i : i \ge 0 \}$. Let
$\alpha(X)$ be an element of End$(M[X])$. The localization map $i: M
\to M_a$ induces a map $i^*:$ End$(M[X]) \to$ End$(M[X]_a) = $
End$(M_a[X])$. We shall denote $i^*(\alpha(X))$ by $\alpha(X)_a$ in
the sequel.}
\end{nt}

\section{\large{Elementary Linear Group $\E_{\varphi}(R)$}}

\begin{de} {\bf Alternating Matrix:} {\rm A matrix from ${\rm M}_n
(R)$ is said to be {\it alternating} if it has the form $\nu - \nu^t$,
where $\nu \in {\rm M}_n(R)$. (It follows that its diagonal elements
are zeros.)}
\end{de}

\begin{de}  {\rm  The group of all invertible $2n \times 2n$ matrices
$\{ \alpha \in \GL_{2n}(R) ~|~ \alpha^t \varphi \alpha = \varphi \}, $
where $\varphi$ is an alternating matrix is called {\it Symplectic
Group} Sp$_{\varphi}(R)$ {\it w.r.t. an invertible alternating matrix
$\varphi$.}  }
\end{de}

\begin{de}  {\rm Let $v \in R^{2n-1}$. Let $\varphi$ be an invertible
alternating matrix of size $2n$ of the form $\big( \begin{smallmatrix}
0 & -c \\ c^t & \nu \end{smallmatrix} \big)$, and $\varphi^{-1}$ be of
the form $\big( \begin{smallmatrix} 0 & d \\ -d^t &
\mu \end{smallmatrix} \big)$, where $c, d \in R^{2n-1}$. In
(\cite{SV}, Lemma 5.4) L.N. Vaserstein considered the matrices
(related to $\varphi$ and $v\in R^{2n-1}$):
\begin{eqnarray*} \alpha &:=& \alpha_\varphi(v) ~ := ~ I_{2n-1} + d^t
v \nu,\\ \beta &:=& \beta_\varphi(v) ~ := ~ I_{2n-1} + \mu v^t c.
\end{eqnarray*}

Note that $\alpha_{\varphi}(v), \beta_{\varphi}(v) \in \E_{2n-1}(R)$
via Lemma \ref{suslin}.

\vskip0.15in

From these matrices he constructed in (\cite{SV}, Lemma 5.4) 
\begin{eqnarray*} C_{\varphi}(v) & = & \begin{pmatrix} 1 & 0 \\ \alpha
v^t & \alpha
\end{pmatrix} ~ = ~ \begin{pmatrix} 1 & 0 \\ v^t &
\alpha \end{pmatrix}  ~~\mbox{and}~~\\ R_{\varphi} (v) & =
& \begin{pmatrix} 1 & v \\ 0 & \beta
\end{pmatrix}.
\end{eqnarray*}

(The notation $C_{\varphi}(v)$, $R_{\varphi}(v)$ is due to us.) In
(\cite{SV}, Lemma 5.4) it is mentioned that these matrices belong to
$\Sp_{\varphi}(R)$. }
\end{de}

\begin{re} {\rm Using this construction, L.N. Vaserstein showed in
(\cite{SV}, Lemma 5.5) that given $\varepsilon \in \E_{2n}(R)$ there
exists $\rho \in \E_{2n-1}(R)$, which is a product of elements of the
type $\alpha_{\varphi}(v), \beta_{\varphi}(v)$, such that $(1 \perp
\rho) \varepsilon \in \Sp_{\varphi}(R)$ (We refer the reader to Lemma
\ref{vStein} for details.)  }
\end{re}

\begin{de} {\bf Elementary Linear Group $\E_{\varphi}(R)$,
$\E_{\varphi}(R,I)$:} {\rm  For an invertible alternating matrix
$\varphi$ of size $2n$ and for $v \in R^{2n-1}$, we have
$\alpha_{\varphi}(v), \beta_{\varphi}(v) \in \E_{2n-1}(R)$ (see Lemma
\ref{suslin}). The subgroup of $\E_{2n-1}(R)$ generated by
$\alpha_{\varphi}(v), \beta_{\varphi}(v)$, for all $v \in R^{2n-1}$ is
called a {\it linear elementary group w.r.t. to an alternating matrix}
$\varphi$ and is denoted by $\E_{\varphi}(R)$.

Let $I$ be an ideal of $R$. The subgroup of $\E_{\varphi}(R)$
generated by $\alpha_{\varphi}(v), \beta_{\varphi}(v)$, for all $v \in
I^{2n-1} (\subseteq R^{2n-1})$ is denoted by $\E_{\varphi}(I)$. And
the normal closure of $\E_{\varphi}(I)$ in $\E_{\varphi}(R)$ is
denoted by $\E_{\varphi}(R,I)$.  }
\end{de}

\begin{lem}  \label{splitting} {\rm (Splitting property)} {\rm For two
row vectors $v, w \in R^{2n-1}$
\begin{eqnarray*} \alpha_{\varphi}(v+w) &=& \alpha_{\varphi}(v)
\alpha_{\varphi}(w) , \\ \beta_{\varphi}(v+w) &=& \beta_{\varphi}(v)
\beta_{\varphi}(w).
\end{eqnarray*} }
Hence the generators of $\E_{\varphi}(R)$ satisfy the splitting
property.
\end{lem}

Proof: Note that
\begin{eqnarray*} \alpha_{\varphi}(v) \alpha_{\varphi}(w) &=&
(I_{2n-1} + d^t v \nu) (I_{2n-1} + d^t w \nu) \\ &=& I_{2n-1} + d^t v
\nu + d^t w \nu ,  ~~~~~~\mbox{{\rm since}  $\nu d^t = 0$}\\  &=&
I_{2n-1} + d^t (v+w) \nu ~=~ \alpha_{\varphi}(v+w),
\end{eqnarray*} 
\begin{eqnarray*} \beta_{\varphi}(v) \beta_{\varphi}(w) &=& (I_{2n-1}
+ \mu v^t c) (I_{2n-1} + \mu w^t c) \\ &=& I_{2n-1} + \mu v^t c + \mu
w^t c,  ~~~~~~\mbox{{\rm since}  $c \mu = 0$}\\ &=&  I_{2n-1} + \mu
(v+w)^t c ~=~ \beta_{\varphi}(v+w).
\end{eqnarray*}
  \hfill{$\square$}

\vskip0.15in

The following lemma is well known.

\begin{lem} \label{conjugation} Let $G$ be a group, and $a_i, b_i \in
G$, for $i = 1, \ldots, n$. Then $\prod_{i=1} ^n a_i b_i = \prod_{i=1}
^n r_i b_i r_i ^{-1} \prod_{i=1} ^n a_i$, where $r_i = \prod_{j=1} ^i
a_j$.  \hfill{$\square$}
\end{lem}

\begin{lem} \label{relative-equal} Let $R[X]$ be the polynomial ring
and $(X)$ be the ideal generated by $X$. Let $\varphi$ be an
invertible alternating matrix of size $2n$. Then $\E_{\varphi}(R[X],
(X)) = \E_{\varphi}(R[X]) \cap \GL_{2n-1}(R[X], (X))$.
\end{lem}

Proof: It is easy to see that $\E_{\varphi}(R[X], (X)) \subseteq
\E_{\varphi}(R[X]) \cap \GL_{2n-1}(R[X], (X))$. The other way
inclusion follows from Lemma \ref{splitting} and Lemma
\ref{conjugation}.  \hfill{$\square$}

\begin{lem} \label{varphi=psi} If $\varphi = \psi_n$, then
$\E_{\varphi}(R) = \E_{2n-1}(R)$.
\end{lem}

Proof: $\E_{\varphi}(R) \subseteq  \E_{2n-1}(R)$ (see Lemma
\ref{suslin}).  Let $e_i$ denote a row vector of length $2n-1$ which
has 1 at the $i^{th}$ place and zeros elsewhere. When $\varphi =
\psi_n$, then for an element $a \in R$ we have
\begin{eqnarray*} E_{12}(a) &=& \alpha_{\psi_n}(a. e_3), \\ E_{1j}(a)
&=&
\begin{cases} \alpha_{\psi_n}(a . e_{j+1}) & \mbox{{\rm if}  $j \ge 3$
{\rm and} $j$ {\rm is even},} \\ \alpha_{\psi_n}(-a . e_{j-1}) &
\mbox{{\rm if} $j \ge 3$ {\rm and} $j$ {\rm is odd},} \\
\end{cases}\\ E_{21}(a) &=& \beta_{\psi_n}(a. e_3), \\ E_{i1}(a)  &=&
\begin{cases} \beta_{\psi_n}(a . e_{i+1}) & \mbox{{\rm if}  $i \ge 3$
{\rm and} $i$ {\rm is even},} \\ \beta_{\psi_n}(- a . e_{i-1}) &
\mbox{{\rm if} $i \ge 3$ {\rm and} $i$ {\rm is odd},} 
\end{cases}
\end{eqnarray*} and hence $\E_{2n-1}(R) \subseteq \E_{\varphi}(R)$.
\hfill{$\square$}

\begin{lem} \label{equality1} Let $\varphi$ and $\varphi^*$ be two
invertible alternating matrices and 
\begin{eqnarray*} \varphi &=& (1 \perp \varepsilon)^t  \varphi^* (1
\perp \varepsilon),
\end{eqnarray*} for some $\varepsilon  \in \E_{2n-1}(R)$. Then
$\E_{\varphi}(R) = \varepsilon^{-1} \E_{\varphi^*}(R) \varepsilon $.
\end{lem}
 
 Proof:  Note that when $\varphi^*$ is of the form
$\big( \begin{smallmatrix} 0 & -c \\ c^t & \nu \end{smallmatrix}
\big)$, and $\varphi{^*{^{-1}}}$ is of the form
$\big( \begin{smallmatrix} 0 & d \\ -d^t & \mu \end{smallmatrix}
\big)$, where $c, d \in R^{2n-1}$, then 
 \begin{eqnarray*} \varphi =  \begin{pmatrix} 0 & -c \varepsilon \\
\varepsilon^t c^t & \varepsilon^t \nu \varepsilon \end{pmatrix}   ~
{\rm and} ~  \varphi^{-1} =  \begin{pmatrix} 0 & d
(\varepsilon^t)^{-1} \\ - \varepsilon^{-1} d^t & \varepsilon^{-1} \mu
(\varepsilon^t)^{-1}  \end{pmatrix} .
 \end{eqnarray*}
  
Using the definition of $\alpha_{\varphi}$ and $\beta_{\varphi}$ we
get
\begin{eqnarray*} \alpha_{\varphi}(v) = Id + (\varepsilon^{-1} d^t) v
(\varepsilon^t \nu \varepsilon)  = \varepsilon^{-1} (Id + d^t (v
\varepsilon^t) \nu) \varepsilon = \varepsilon^{-1}
\alpha_{\varphi^*}(v \varepsilon^t) \varepsilon, \\ \beta_{\varphi}(v)
= Id + (\varepsilon^{-1} \mu (\varepsilon^t)^{-1})  v^t (c
\varepsilon) = \varepsilon^{-1} ( Id + \mu ((\varepsilon^t)^{-1} v^t)
c ) \varepsilon = \varepsilon^{-1} \beta_{\varphi^*}(v
\varepsilon^{-1}) \varepsilon,
\end{eqnarray*} and hence the equality follows.  \hfill{$\square$}

 \begin{lem} \label{equality2} Let $\varphi = (1 \perp \varepsilon)^t
~ \psi_n (1 \perp \varepsilon)$, for some $\varepsilon  \in
\E_{2n-1}(R)$. Then $\E_{\varphi}(R) = \E_{2n-1}(R)$.
 \end{lem}
 
Proof: Follows from the above two lemmas.  \hfill{$\square$}

\section{\large{Dilation and LG Principle for $\E_{\varphi\otimes
R[X]}(R[X])$}}

Here we establish dilation principle and Local-Global principle for
$\E_{\varphi \otimes R[X]}(R[X])$. 

\begin{lem} \label{dilation} {\rm (Dilation principle)} Let $\varphi$
be an alternating matrix of Pfaffian 1 of size $2n$, with $n \ge
2$. Let $a \in R$ be a non-nilpotent element, and let $\varphi = (1
\perp \varepsilon)^t ~ \psi_n ~ (1 \perp \varepsilon)$, for some
$\varepsilon \in \E_{2n-1}(R_a)$ over the ring $R_a$. Let $\theta(X)
\in \E_{\varphi \otimes R_a[X]}(R_a[X])$, with $\theta(0)= Id$. Then
there exists $\theta^*(X) \in \E_{\varphi \otimes R[X]}(R[X])$ such
that $\theta^*(X)$ localises to $\theta(bX)$, for some $b \in (a^d)$,
$d \gg 0$, and $\theta^*(0)=Id$.
\end{lem}

Proof: We are given that $\theta(X) \in \E_{\varphi \otimes
R_a[X]}(R_a[X])$, where $\varphi = (1 \perp \varepsilon)^t ~ \psi_n ~
(1 \perp \varepsilon)$, for some $\varepsilon \in \E_{2n-1}(R_a)$ over
the ring $R_a$. Therefore by Lemma \ref{equality1} we have $\theta(X)
= \varepsilon^{-1} \eta(X) \varepsilon$, for some $\eta(X) \in
\E_{2n-1}(R_a[X])$. Since $\eta(0) = Id.$ we can write $\eta(X) =
\prod \gamma_l ~E_{i_l j_l}(X f_l(X)) ~\gamma_l ^{-1}$, where
$\gamma_l \in \E_{2n-1}(R_a)$, and $f_l(X) \in R_a[X]$ (see Lemma
\ref{conjugation}). Using
commutator identities for the generators of the elementary linear
group we get $\eta(Y^r X) = \prod E_{i_k j_k}(Yh_k(X,Y)/ a^s)$, for
large integer $r$. Here $h_k(X,Y) \in R[X,Y]$ and either $i_k =1$ or
$j_k =1$. Using equalities appearing in the Lemma \ref{varphi=psi} it
is clear that $\eta(Y^r X)$ is product of the elements of the form
$\alpha_{\psi_n}((Yh_k(X,Y)/ a^s). e_i)$ or
$\beta_{\psi_n}((Yh_k(X,Y)/ a^s). e_j)$, where $2 \le i, j \le 2n-1$. 

Note that $\alpha_{\psi_n}((Yh_k(X,Y)/ a^s). e_i) = \varepsilon
\alpha_{\varphi}((Yh_k(X,Y)/ a^s). e_i (\varepsilon^t)^{-1})
\varepsilon^{-1}$ and \\$\beta_{\psi_n}((Yh_k(X,Y)/ a^s). e_j) =
\varepsilon  \beta_{\varphi}((Yh_k(X,Y)/ a^s). e_j \varepsilon)
\varepsilon^{-1}$. Therefore $\theta(Y^rX)$ is product of elements of
the form $\alpha_{\varphi}((Yh_k(X,Y)/ a^s). e_i
(\varepsilon^t)^{-1})$ or $\beta_{\varphi}((Yh_k(X,Y)/ a^s). e_j
\varepsilon)$. Let $t$ be the maximum power of $a$ appearing in the
denominators of $\varepsilon$ and $(\varepsilon^t)^{-1}$. Set $d =
s+t$. Define $\theta^*(X,Y)$ as product of elements of the form
$\alpha_{\varphi}(Yh_k(X, a^dY). a^t e_i (\varepsilon^t)^{-1})$ and
$\beta_{\varphi}(Yh_k(X, a^dY) a^t. e_j \varepsilon)$. Note that
$\theta^*(X, Y) \in \E_{\varphi \otimes R[X, Y]}(R[X, Y])$. We obtain
$\theta^*(X)$ substituting $Y=1$ in $\theta^*(X,Y)$. Clearly
$\theta^*(X)$ localises to $\theta(bX)$ for some $b \in (a^d)$, and
$\theta^*(0) = Id$.  \hfill{$\square$}

\begin{re} \label{local} Let $(R, \gm)$ be a local ring and $\varphi$
be an alternating matrix of Pfaffian 1 over $R$ of size $2n$. Then
$\varphi = \varepsilon^t  \psi_n \varepsilon$, for some $\varepsilon
\in \E_{2n}(R)$.
 \end{re}
  
 We recollect an observation of Rao-Swan stated in the introduction of
\cite{SwR}. We make a contextual observation which the proof shows
and include it for completeness. 
 
 \begin{lem} \label{vStein} {\rm (Rao-Swan)} Let $n \ge 2$ and
$\varepsilon \in \E_{2n}(R)$. Then there exists $\rho \in
\E_{\psi_n}(R) \subseteq \E_{2n-1}(R)$ such that $(1 \perp \rho)
\varepsilon \in \ESp_{2n}(R)$.
\end{lem}

Proof: Let $\varepsilon = \varepsilon_r \ldots \varepsilon_1$, where
each $\varepsilon_i$ is of the form $ \big( \begin{smallmatrix} 1 &
v\\ 0 & I \end{smallmatrix}  \big)$ or $ \big( \begin{smallmatrix} 1 &
0\\ v^t & I \end{smallmatrix} \big)$, where $v=(a_1, \ldots, a_{2n-1})
\in R^{2n-1}$ (see \cite{SV}, Lemma 2.7).  We prove the result using
induction on $r$. It is clear when $r=0$. Let $r \ge 1$. Let us assume
the result is true for $r-1$, i.e, when $\varepsilon =
\varepsilon_{r-1} \ldots \varepsilon_{1}$, then there exists a $\delta
\in \E_{2n-1}(R)$ such that $(1 \perp \delta) \varepsilon \in
\ESp_{2n}(R)$. We will prove the result when number of generators of
$\varepsilon$ is $r$. We have
\begin{eqnarray*} C_{\psi_n}(v) & = &
\begin{pmatrix} 1 & 0\\ 0 & \alpha \end{pmatrix} \begin{pmatrix} 1 & 0
\\  v^t & I_{2n-1} \end{pmatrix} ~ = ~ \prod_{i=2}^{2n}
se_{i1}(a_{i-1}),\\ R_{\psi_n}(v) & = &
\begin{pmatrix} 1 & 0 \\ 0 & \beta \end{pmatrix}  \begin{pmatrix} 1 &
v \\  0 & I_{2n-1} \end{pmatrix} ~ = ~ \prod_{i=2}^{2n}
se_{1i}(a_{i-1}).
\end{eqnarray*}

Note that $\alpha= \alpha_{\psi_n}(v), \beta=\beta_{\psi_n}(v) \in
\E_{2n-1}(R)$. Let us set $\gamma$ equal to either $\alpha$ or $\beta$
depending on the form of $\varepsilon_1$. Now,
$\big( \begin{smallmatrix} 1 & 0 \\ 0 & \gamma \end{smallmatrix}
\big) \varepsilon_1  \in \ESp_{2n}(R)$, and each $\eta_i =
\big( \begin{smallmatrix} 1 & 0 \\ 0 & \gamma \end{smallmatrix} \big)
\varepsilon_i  \big( \begin{smallmatrix} 1 & 0 \\ 0 &
\gamma^{-1} \end{smallmatrix} \big)$ is of the form
$\big( \begin{smallmatrix} 1 & v\\ 0 & I \end{smallmatrix} \big)$ or
$\big( \begin{smallmatrix} 1 & 0\\ v^t & I \end{smallmatrix}
\big)$. Now we have 
\begin{eqnarray*} \varepsilon &=& \begin{pmatrix} 1 & 0 \\ 0 &
\gamma^{-1} \end{pmatrix} \eta_r \ldots \eta_2 \begin{pmatrix} 1 & 0
\\ 0 & \gamma 
\end{pmatrix}\varepsilon_1 .
\end{eqnarray*}

By induction hypothesis $(1 \perp \delta) \eta_r \ldots \eta_2 \in
\ESp_{2n}(R)$, for some $\delta \in \E_{2n-1}(R)$. Hence $(1 \perp
\rho) \varepsilon \in \ESp_{2n}(R)$, where $\rho = \delta^{-1} \gamma
\in \E_{2n-1}(R)$.   \hfill{$\square$}

  \begin{cor} {\rm (Rao-Swan)} \label{rao-swan}  For $n \ge 2$ and 
$\varepsilon \in
\E_{2n}(R)$, we have an $\varepsilon_0 \in \E_{\psi_n}(R) \subseteq
\E_{2n-1}(R)$ such that $\varepsilon^t  \psi_n \varepsilon = (1 \perp
\varepsilon_0)^t  \psi_n (1 \perp \varepsilon_0)$.
 \end{cor}

 Proof: Using Lemma \ref{vStein} we get $\varepsilon_0 \in
\E_{\psi_n}(R) \subseteq \E_{2n-1}(R)$ such that $(1
\perp\varepsilon_0) \varepsilon^{-1} \in \ESp_{2n}(R)$, and hence
${\varepsilon^{-1}}^t (1 \perp \varepsilon_0)^t ~ \psi_n ~ (1
\perp\varepsilon_0) \varepsilon^{-1} = \psi_n$.

Therefore we have
\begin{eqnarray*} \varepsilon^t \psi_n \varepsilon &=& \varepsilon^t ~
\{ {\varepsilon^{-1}}^t  (1 \perp \varepsilon_0)^t \} ~ \psi_n ~ \{ (1
\perp\varepsilon_0)  \varepsilon^{-1} \} ~ \varepsilon\\ &=& (1 \perp
\varepsilon_0)^t \psi_n (1 \perp \varepsilon_0).
\end{eqnarray*} \hfill{$\square$}

\begin{cor} \label{localRing} Let $(R, \gm)$ be a local ring and
$\varphi$ be an alternating matrix of Pfaffian 1 over $R$ of size $2n$
with $n \ge 2$. Then the groups $\E_{2n-1}(R)$ and $\E_{\varphi}(R)$
are equal.
\end{cor}

Proof: Over the local ring $(R, \gm)$ we have $\varphi = (1 \perp
\varepsilon_0)^t \psi_n (1 \perp \varepsilon_0)$, for some
$\varepsilon_0 \in \E_{2n-1}(R)$. This follows by Remark \ref{local}
and Corollary \ref{rao-swan}. Using Lemma \ref{equality1} we get
$\E_{\varphi}(R) = \varepsilon_0^{-1} \E_{\psi_n}(R)
\varepsilon_0$. Using Lemma \ref{varphi=psi} we get $\E_{\psi_n}(R) =
\E_{2n-1}(R)$, and hence the equality follows.  \hfill{$\square$}

\vskip0.15in

Using dilation principle we prove the following variant of
D. Quillen's Local-Global principle (see \cite{Q}). The argument is
standard. We include the proof for completeness. Interested reader may
refer to \cite{bkr} for a survey about the Local-Global principle.

\begin{thm} \label{LG} {\rm (Local-Global principle)} Let $\varphi$ be
an alternating matrix of Pfaffian $1$ of size $2n$, with $n \ge
2$. Let $\theta(X) \in \SL_{2n-1}(R[X])$, with $\theta(0) = Id$. If
$\theta(X)_\gm \in \E_{\varphi \otimes R_\gm[X]}(R_\gm[X])$, for all
maximal ideal $\gm$ of $R$, then $\theta(X) \in \E_{\varphi \otimes
R[X]}(R[X])$.
\end{thm}

Proof:   For each maximal ideal $\gm$ of $R$ one can suitably choose
an element $a_\gm$ from $R \setminus \gm$ such that $\theta(X)_{a_\gm}
\in \E_{\varphi \otimes R_{a_\gm}[X]}(R_{a_\gm}[X])$ and also one has
$\varphi = (1 \perp \varepsilon)^t ~ \psi_n ~ (1 \perp \varepsilon)$,
for some $\varepsilon \in \E_{2n-1}(R_{a_\gm})$. Define $\gamma(X,Y)
=\theta(X+Y)_{a_\gm} \theta(Y)_{a_\gm}^{-1}$.  It is clear that
\begin{eqnarray*}  \gamma(X,Y) &\in& \E_{\varphi \otimes
R_{a_\gm}[X,Y]}(R_{a_\gm}[X,Y])
\end{eqnarray*} and $\gamma(0,Y) = Id$. Therefore $\gamma(b_\gm X,Y)
\in \E_{\varphi \otimes R[X,Y]}(R[X,Y])$, where $b_\gm \in (a_m^d)$
for $d \gg 0$ (see Lemma \ref{dilation}). Note that the ideal
generated by $a_\gm^d$'s is the whole ring $R$. Therefore, $c_1
a_{\gm_1} ^d+ \cdots + c_k a_{\gm_k} ^d = 1 $, for some $c_i \in
R$. Let $b_{m_i} = c_i a_{m_i} ^d \in (a_{m_i} ^d)$. It is easy to see
that $\theta(X) = \prod_{i=1} ^{k-1} \gamma(b_{m_i}X, T_i)
\gamma(b_{m_k}X, 0)$, where $T_i = b_{m_{i+1}}X + \cdots +
b_{m_k}X$. Each term in the right hand side of this expression belongs
to $\E_{\varphi \otimes R[X]}(R[X])$, and hence $\theta(X) \in
\E_{\varphi \otimes R[X]}(R[X])$.  \hfill{$\square$}

\vskip0.15in

We recall Swan-Weibel's trick to establish the Local-Global principle
in the graded case.

\begin{thm} \label{LG-multi} {\rm (Graded case of Local-Global
principle)} Let $\varphi$ be an alternating matrix of Pfaffian $1$ of
size $2n$, with $n \ge 2$. Let $\theta(X_1, \ldots, X_t) \in
\SL_{2n-1}(R[X_1, \ldots, X_t])$, with $\theta(0, \ldots, 0) = Id$. If
$\theta(X_1, \ldots, X_t)_\gm \in \E_{\varphi \otimes R_\gm[X_1,
\ldots, X_t]}(R_\gm[X_1, \ldots, X_t])$, for all maximal ideal $\gm$
of $R$, then $\theta(X_1, \ldots, X_t) \in \E_{\varphi \otimes R[X_1,
\ldots, X_t]}(R[X_1, \ldots, X_t])$.
\end{thm}

Proof: Let us denote $S = R[X_1, \ldots, X_t]$. Note that $S$ is a
graded ring with the grading $S = S_0 \oplus S_1 \oplus S_2 \oplus
\cdots$, and $S_0 = R$. Consider the ring homomorphism $f : S \lra
S[T]$ given by $f(a_0 + a_1 + a_2 + \cdots) = a_0 + a_1 T + a_2 T^2 +
\cdots$, where each $a_i$ is a homogeneous component belongs to
$S_i$. Let us denote $\theta(X_1, \ldots, X_t) = (\theta_{i j})$,
where $\theta_{i j} \in S$. We set $\tilde{\theta}(T) = ( f(\theta_{i
j}))$. Note that $\tilde{\theta}(1) = (\theta_{i j})$, and
$\tilde{\theta}(0) = \theta(0, \ldots, 0) = Id$.

Let $\gm_0$ be a maximal ideal of $R$ and let $M_0 = R \setminus
\gm_0$. Since $(\theta_{ij})_{M_0} \in \E_{\varphi}(S_{M_0})$, we have
$\tilde{\theta}(T)_{M_0} \in \E_{\varphi \otimes
S_{M_0}[T]}(S_{M_0}[T])$. Therefore, there is a $s_{m_0} \in M_0$ such
that $\tilde{\theta}(T)_{s_{m_0}} \in \E_{\varphi \otimes
S_{s_{m_0}}[T]}(S_{s_{m_0}}[T])$. If $\gm$ is a maximal ideal of $S$
then $s_{m_0} \notin \gm$ for some $\gm_0$. Therefore,
$\tilde{\theta}(T)_\gm \in \E_{\varphi \otimes S_\gm[T]}(S_\gm[T])$,
for all maximal ideals $\gm$ of $S$. Moreover, the ideal generated by
all $s_{m_0}$, for all maximal ideals $\gm_0$ of $R$, is the whole
ring $R$. Hence, $\tilde{\theta}(T) \in \E_{\varphi}(S[T])$. Now
substituting $T=1$ we get $\tilde{\theta}(1) = (\theta_{i j}) =
\theta(X_1, \ldots, X_t) \in \E_{\varphi \otimes R[X_1, \ldots,
X_t]}(R[X_1, \ldots, X_t])$.  \hfill{$\square$}

\section{\large{Equality of $\E_{2n-1}(R)$ and $\E_{\varphi}(R)$}}

We are now ready to prove the main theorem of this note.
\vskip0.15in

\begin{thm} \label{equalityGroup} Let $\varphi$ be an alternating
matrix of Pfaffian $1$ of size $2n$, with $n \ge 2$. Then the groups
$\E_{2n-1}(R)$ and $\E_{\varphi}(R)$ are equal.
\end{thm}

Proof: By Lemma \ref{suslin}, it follows that $\E_{\varphi}(R)
\subseteq \E_{2n-1}(R)$. Let $\lambda \in \E_{2n-1}(R)$. Then there
exists $\lambda(X) \in \E_{2n-1}(R[X])$ such that $\lambda(1) =
\lambda$ and $\lambda(0) = Id$. 
 
Let $\gm$ be a maximal ideal of $R$ and  $R_\gm$ be the local ring at
$\gm$. Over the local ring $R_\gm$, we have $\varphi = (1 \perp
\varepsilon)^t \psi_n (1 \perp \varepsilon)$, for some $\varepsilon
\in \E_{2n-1}(R_\gm)$. Therefore, for each maximal ideal $\gm$ of $R$,
we have $\lambda(X)_\gm \in \E_{2n-1}(R_\gm[X]) = \E_{\varphi \otimes
R_\gm[X]}(R_\gm[X])$ (see Corollary \ref{localRing}). By Theorem
\ref{LG} it follows that $\lambda(X) \in \E_{\varphi \otimes
R[X]}(R[X])$. Put $X=1$ to get $\lambda = \lambda(1) \in
\E_{\varphi}(R)$, i.e, $\E_{2n-1}(R) \subseteq \E_{\varphi}(R)$. Hence
the equality follows.  \hfill{$\square$}

\vskip0.15in

Using the ideas of the proof of (\cite{rao}, Lemma 3.1) we can prove a
relative version of the  above theorem w.r.t. an ideal.

\begin{thm} \label{equalityGroup-rel} Let $I$ be an ideal of the ring
$R$. Let $\varphi$ be an alternating matrix of Pfaffian $1$ of size
$2n$, with $n \ge 2$. Then the groups $\E_{2n-1}(R, I)$ and
$\E_{\varphi}(R, I)$ are equal.
\end{thm}

Proof: It is clear that $\E_{\varphi}(R, I) \subseteq \E_{2n-1}(R, I)$
(see Lemma \ref{suslin}). Let $\lambda \in \E_{2n-1}(R, I)$ and
$\lambda = \prod_{k=1} ^t \gamma_k E_{i_k j_k}(a_k) \gamma_k^{-1}$,
where $\gamma_k \in \E_{2n-1}(R)$ and $a_k \in I$. We define
$\Lambda(X_1 , \ldots, X_t) = \prod_{k=1} ^t \gamma_k E_{i_k j_k}(X_k)
\gamma_k^{-1}$, which is in $\E_{2n-1}(R[X_1, \ldots, X_t], (X_1,
\ldots, X_t))$.

By Corollary \ref{localRing} for all maximal ideals $\gm$ of $R$ we
have 
\begin{eqnarray*} \E_{2n-1}(R_\gm[X_1, \ldots, X_t]) &=& \E_{\varphi
\otimes R_\gm[X_1, \ldots, X_t]}(R_\gm[X_1, \ldots, X_t]).
\end{eqnarray*}  Therefore, $\Lambda(X_1, \ldots, X_t)_\gm \in
\E_{\varphi \otimes R_\gm[X_1, \ldots, X_t]}(R_\gm[X_1, \ldots,
X_t])$, for all maximal ideals $\gm$ of $R$. Hence, $\Lambda(X_1,
\ldots, X_t) \in \E_{\varphi \otimes R[X_1, \ldots, X_t]}(R[X_1,
\ldots, X_t])$ by Theorem \ref{LG-multi}. Also, $\Lambda(X_1, \ldots,
X_t) \in \SL_{2n-1}(R[X_1, \ldots, X_t], (X_1, \ldots, X_t))$ and
hence $\Lambda(X_1, \ldots, X_t)$ is in the relative group
$\E_{\varphi \otimes R[X_1, \ldots, X_t]}(R[X_1, \ldots, X_t], (X_1,
\ldots, X_t))$ (see Lemma \ref{relative-equal}). Substituting $(X_1,
\ldots, X_t) = (a_1, \ldots, a_t)$ we get $\lambda \in \E_{\varphi}(R,
I)$.  \hfill{$\square$}

 \begin{re} \label{alternating} {\rm The condition that the
alternating matrices, in this article, are of Pfaffian one can be
extended to all invertible alternating matrices by observing that an
invertible alternating matrix over a local ring which is congruent to
$(u~\psi_1~\perp~\psi_{n-1}) \\~({\rm mod}~I)$, where $u$ = Pfaffian
$\varphi$, is of the form $(1 \perp E)^t (u \psi_1 \perp \psi_{n-1})
(1 \perp E)$, for some relative elementary matrix $E$. }
\end{re}

\vskip0.15in

\noindent {\bf Acknowledgement:} The first named author thanks Department of 
Science and Technology, Govt. of India for awarding INSPIRE Faculty Fellowship 
which supported this work. She also thanks Tata Institute of Fundamental 
Research, Mumbai for its hospitality.

\end{document}